\documentclass[11pt]{article}
\usepackage[utf8]{inputenc}
\usepackage{tabularx,lscape,longtable}
\usepackage[english]{babel}
\usepackage[dvipsnames]{xcolor}
\usepackage{amsfonts,comment}
\usepackage{multicol}
\usepackage{amssymb,latexsym,amsmath}
\makeatletter \oddsidemargin 0in \evensidemargin 0in \textwidth
16cm \RequirePackage[dvips]{graphicx} \textheight 20cm
\DeclareGraphicsExtensions{.jpg, .eps}
\newtheorem{p1}{Proposition}[section]
\newtheorem{l1}{Lemma}[section]
\newtheorem{c1}{Corollary}[section]
\newtheorem{d1}{Definition}[section]
\newtheorem{r1}{Remark}[section]

\definecolor{cadmiumred}{rgb}{0.89,0.0,0.13}
\renewcommand{\baselinestretch}{1.3}
\renewcommand{\theequation}{\thesection.\arabic{equation}}
\begin{document}
\title{\bf
Use of stochastic orders and statistical dependence in error analysis for multi-component system}
{\author{Subarna Bhattacharjee$^{1}$$^{}$\thanks{Corresponding author~:
E-mail ID: subarna.bhatt@gmail.com}, Aninda Kumar Nanda$^{2}$, Subhashree Patra$^{3}$ \\
{\it $^{1,3}$ Department of Mathematics,Ravenshaw University, Cuttack-753003, Odisha, India}\\
{\it $^{2}$ Indian Statistical Institute, Delhi Center. India}\\}
\date{\today}
\maketitle
\begin{abstract}
In this paper, we analyze the relative errors that crop up in the various reliability measures due to the tacit assumption that the components are independently working associated with a $n$-component series system or a parallel system where the components are dependent and follow a well- defined multivariate Weibull or exponential distribution. We also list some important observations which the previous authors have not noted in their earlier works. In this paper, we focus on the incurred error in multi-component series and parallel systems having multivariate Weibull distributions. In the upcoming sections, we establish that the present study has relevance with stochastic orders and statistical dependence which were not previously pointed out by previous authors. 
\end{abstract}
\section{Introduction}
${}$\hspace{0.5cm}To describe the aging phenomena of a system or a unit with lifetime represented by a non-negative random variable $T,$ we commonly use the functions, namely survival function (SF), $\overline{F}(t)=1-F(t);$ failure rate function ($FR$), $r(t)=f(t)/\overline{F}(t);$ reversed failure rate function ($RFR$), $\mu(t)=f(t)/F(t);$ mean residual life function ($MRL$), $m(t)=\frac{\int_{t}^{\infty}\overline{F}(x)dx}{\overline{F}(t)}$ defined at time $t\geq 0.$ Here, $f(t)$ and $F(t)$ are respectively the density function and distribution functions of $T$.
Recently, aging intensity function ($AI$) (cf. Jiang et al. (2003), Nanda et al. (2007), Bhattacharjee (2013) etc.) widely used in aging analysis is defined as
$$L(t)=\frac{-tf(t)}{\overline{F}(t)\ln \overline{F}(t)}=\frac{tr(t)}{\int_{0}^{t}r(u)du}, t>0.$$
\hspace*{1cm} A system with a well defined statistical distribution is said to belong to the non-parametric aging class of increasing (decreasing) failure rate, increasing (decreasing) failure rate average denoted by IFR (DFR) and IFRA (DFRA) according as $r(t)$ or $\frac{1}{t}\int_{0}^{t}r(x)dx$ is increasing (decreasing) in $t>0.$ Other aging classes having increasing (decreasing) aging intensity function called as increasing (decreasing) aging intensity classes also  exist  in literature. The corresponding aging classes are denoted, respectively, by $IAI$ ($DAI$).\\
\hspace*{1cm}In reliability theory, we come across different $n$ component structures such as series, parallel and $k$ out-of-$n$ systems. During the study of the reliability measures of any $n$-component system, we most often assume that all the components are independently working. Such an assumption eases mathematical complexity while analyzing the system behavior, and hence is very common. However, in
reality, the lifetimes of the components are very interdependent. Such an assumption of
independence hence leads to inaccurate analysis of data.
The analysis of incurred error in computation of various reliability functions arising from
a faulty assumption of components being independently working has been given importance
by various authors in recent past. Authors viz.,  Moeschberger and Klein (1984), Klein and Moeschberger (1986, 1987), Gupta and Gupta (1990), and Nanda et al.(2022) focused on
the relative errors of some bivariate distributions. We now give a brief outline on the works done by aforementioned authors. \\
\hspace*{1cm}Moeschberger and Klein (1984) have studied the consequences of deviations from independence when the lifetimes of the components of a series system are exponentially distributed. They have assumed the joint distribution to follow a Gumbel bivariate exponential model.\\
\hspace*{1cm} Klein and Moeschberger (1986) studied the magnitude of the errors under a similar assumption about the component lifetimes to have independent exponential distributions when, in fact, the lifetimes follow the bivariate exponential distribution of Marshall and Olkin (series or parallel systems) or that of Freund (parallel systems).\\
\hspace*{1cm} Klein and Moeschberger (1987) investigated on a similar line by assuming the joint distribution to follow either one of the three Gumbel bivariate exponential models, the Downton bivariate exponential model, or the Oakes bivariate exponential model. They also concluded that the amount of error incurred depends upon the correlation between lifetimes and the relative mean life of the two components.\\
\hspace*{1cm} Gupta and Gupta (1990) studied relative error in reliability measures such as the reliability function, the failure rate and the mean residual life under the erroneous assumption of independence when in fact lifetimes follow a bivariate exponential model. The behaviour of these errors has been discussed to examine their structure as a function of time.\\
\hspace*{1cm} Recently, Nanda et al. (2022) focused on relative errors in the various reliability measures incurred due to the assumption of independence of components in a series of two components when components actually follow bivariate exponential distributions. They studied relative errors in four reliability measures, viz. Reliability function, failure rate function, mean residual life function, and reversed failure rate function.\\
${}$\hspace{0.5cm}We now give a sketch of the  organization of this paper. In Section 2, we highlight the importance of statistical dependence and stochastic orders in the study of error analysis incurred due to deviation from dependence assumption of various reliability functions for multicomponent systems. In particular, we place our specific findings related to two-component series and parallel systems. In Section 3, we review error assessment in the $n$-component series system when components follow different multivariate exponential distributions but components are mistakenly assumed to be independently working.  Subsequently, we study a  generalized version, namely when components follow some multivariate Weibull distributions. \\  \hspace*{1cm}This study finally proves that the nature of error assessment becomes easily tractable if we use this alternative approach involving stochastic orders and statistical dependence which earlier authors have either not noted or  not pointed out in their works.
\section{Role of statistical dependence and stochastic orders in error analysis of reliability functions}
 To avoid such error in results, instead of independent distribution, many authors have preferred a number of distributions like, Marshall and Olkin’s multivariate exponential, Multivariate Gumbel’s type 1, Marshall and Olkin’s multivariate Weibull (1967a), multivariate Crowder (1989), multivariate Lee (1979), Lu and Bhattacharyya (1990) I, Farlie-Gumbel-Morgenstern Weibull, Lu and Bhattacharyya (1990) II, and many others. However, if we have some knowledge on the behavior of errors or some estimate on the error bound, we could decide if we assume independence and prefer mathematical simplicity (if we know the error is within our allowed limit), or keep the mathematical complexity and get accurate results without assuming independence.\\
 \hspace*{1cm}  We know that errors occur in different functions, viz. $SF$, $FR$, $RFR$, $MRL$, and $AI$ due to the faulty assumption that independent components comprise a system, where they are actually dependent. In this paper, we denote the corresponding relative errors in the aforementioned functions respectively by $E^{\overline{F}}(t)$, $E^{r}(t)$,  $E^{\mu}(t)$,  $E^{L}(t)$ at time $t\geq 0.$ Clearly, for a series system with lifetime $T^{S}$, we evaluate and denote the corresponding relative errors as $$E_{S}^{\overline{F}}(t)=\frac{\overline{F}_{D}^{S}(t)-\overline{F}_{I}^{S}(t)}{\overline{F}_{I}^{S}(t)}, E_{S}^{r}(t)=\frac{r_{D}^{S}(t)-r_{I}^{S}(t)}{r_{I}(t)}, E_{S}^{\mu}(t)=\frac{\mu_{D}^{S}(t)-\mu_{I}^{S}(t)}{\mu_{I}(t)}, E_{S}^{L}(t)=\frac{L_{D}^{S}(t)-L_{I}^{S}(t)}{L_{I}(t)},$$ where \begin{equation}
\label{seriesd}
    \overline{F}_{D}^{S}(t)=P(T_{D}^{S} >t)=P(T_1>t,T_2>t,\ldots,T_n>t)
\end{equation} and 
\begin{equation}
\label{seriesi}
    \overline{F}_{I}^{S}(t)=P(T_{I}^{S} >t)=\prod_{i=1}^{n}P(T_i>t).
\end{equation}
 As noted, we have denoted the relative errors in $SF, FR, RFR, AI$ of the series system by $E^{\overline{F}}_{S}(t)$, $E^{r}_{S}(t)$, $E^{\mu}_{S}(t)$, $E^{L}_{S}(t)$, respectively, at time $t> 0.$. The reader will understand that in the evaluation of the relative error in $SF$ ($E_{S}^{\overline{F}}(t)$) of the series system with lifetime $T_{D}^{S}$, $SF$ of the series system ( $\overline{F}_{D}^{S}(t)$ ) with dependent components and $SF$  of the same series system with lifetime $T_{I}^{S}$ ($\overline{F}_{I}^{S}(t)$) with independent components are taken into account. The same analogy holds for other reliability functions.\\
\hspace*{1cm}In case, we are dealing with a parallel system, the superscript $S$ will be replaced by $P$ wherever used, and as a result, the functions considered correspond to a parallel system. Similarly, for a parallel system with lifetime denoted by $T^{P}$, we evaluate  the corresponding relative 
errors as 
\begin{eqnarray}
E_{P}^{\overline{F}}(t)&=&\frac{\overline{F}_{D}^{P}(t)-\overline{F}_{I}^{P}(t)}{\overline{F}_{I}^{P}(t)},\nonumber\\
&=&\frac{{F}_{I}^{P}(t)-{F}_{D}^{P}(t)}{\overline{F}_{I}^{P}(t)}
\end{eqnarray}$$E_{P}^{r}(t)=\frac{r_{D}^{P}(t)-r_{I}^{P}(t)}{r_{I}(t)}, E_{P}^{\mu}(t)=\frac{\mu_{D}^{P}(t)-\mu_{I}^{P}(t)}{\mu_{I}(t)}, E_{P}^{L}(t)=\frac{L_{D}^{P}(t)-L_{I}^{P}(t)}{L_{I}(t)},$$ where $$\overline{F}_{D}^{P}(t)=1-P(T^{P}_{D}\leq t)=1-P(T_1\leq t,T_2\leq t,\ldots,T_n\leq t)$$ and $$\overline{F}_{I}^{P}(t)=1-P(T^{P}_{I}\leq t)=1-\prod_{i=1}^{n}P(T_i\leq t),$$ where the symbols used have direct analogy as that used for series system.\\
If the relative error ($E$) corresponding to any reliability function, say $SF$, $FR$, $RFR$ etc. is positive (negative), then we say there is under-assessment (over-assessment), denoted by $UA$ ($OA$) arising due to faulty assumption of independence of components. In particular, if $E_{S}^{\overline{F}}(t) \geq (\leq) 0,$ then there is $UA$ ($OA$) in survival function. The same convention is extended for other functions, in discussion.\\
\hspace*{1cm}To infer about $OA$ or $UA,$ it is sufficient to focus our attention on the numerator of each of the formulas involving relative error, since the denominator involves a non-negative function.
\subsection{Use of Statistical dependence in error analysis}
Now, we show how error analysis, as discussed in the text, is related to the study of statistical dependence. For this, 
we now recall a couple of definitions from Nelson  (2006) and other references therein.
\begin{d1}\label{r1} Let $T_1,T_2,\ldots,T_n$ be random variables. $T_1,T_2,\ldots,T_n$ are said to be positively (negatively) upper orthant dependent
PUOD (NUOD) if $$P(T_{1}>t_{1}, T_{2}>t_{2},\ldots,T_{n}>t_{n}) \geq (\leq) \prod_{i=1}^{n}P(T_{i}>t_{i}), $$ i.e., $\overline{F}_{T_1,T_2,\ldots, T_n} (t_1,t_2,\ldots,t_n)\geq (\leq) \prod_{i=1}^{n}\overline{F}_{T_{i}}(t_i),$ where $(t_1,t_2,\ldots,t_n)\in R^{n}.$ \\Clearly, $\overline{F}_{T_1,T_2,\ldots, T_n}(\cdot,\cdot,\ldots,\cdot)$ , $\overline{F}_{T_{i}}(\cdot)$ represent the joint and marginal survival  function of random variables under consideration.
\end{d1}
\begin{d1}\label{r2}
Let $T_1,T_2,\ldots,T_n$ be random variables. $T_1,T_2,\ldots,T_n$ are said to be  positively (negatively) lower orthant dependent
PLOD (NLOD) if $$P(T_{1}\leq t_{1}, T_{2}\leq t_{2},\ldots,T_{n}\leq t_{n}) \geq (\leq) \prod_{i=1}^{n}P(T_{i}\leq t_{i}), $$ i.e., ${F}_{T_1,T_2,\ldots, T_n} (t_1,t_2,\ldots,t_n)\geq (\leq) \prod_{i=1}^{n}{F}_{T_{i}}(t_i)$ where $(t_1,t_2,\ldots,t_n)\in R^{n}.$
\end{d1}
The next proposition immediately reveals how statistical dependence is related to $OA$ and $UA$. An outline of the proof is given for immediate follow-up.
\begin{p1}
Let $T_1, T_2,\ldots,T_n$ denote the component lifetimes of a series (parallel) system and let relative error in $SF$ due to faulty assumption of components independently working in the series (parallel) system be denoted by $E^{\overline{F}}_{S}(t)$ ($E^{\overline{F}}_{P}(t)$). Then \begin{enumerate}
    \item[$(i)$]  $OA (UA)$ in error assessment of $SF$ for an $n$ component series system is equivalent to the fact that $T_1, T_2,\ldots,T_n$ are PUOD (NUOD).
\item [$(ii)$]  $OA (UA)$ in error assessment of $SF$ for an $n$ component parallel system is equivalent to the fact that $T_1, T_2,\ldots,T_n$ are PLOD (NLOD).
\end{enumerate}
\end{p1}
{\bf Proof.} Clearly, $E^{\overline{F}}_{S}(t)\geq (\leq)~ 0,$ if and only if
$\overline{F}_{D}^{S}(t)\geq (\leq)  
 \overline{F}_{I}^{S}(t)$ for all $t\geq 0.$ Similarly, $E^{\overline{F}}_{P}(t)\geq (\leq)~ 0,$ if and only if
$\overline{F}_{D}^{P}(t)\geq (\leq)  \overline{F}_{I}^{P}(t)$ for all $t\geq 0.$  The proof is easy to follow from Definition \ref{r1} and \ref{r2}. $\hfill\Box$\\
We know that a bivariate distribution bears the same inequality (in particular, greater inequality) between the joint distribution (joint survival function) and that of the product of marginal distributions (marginal survival functions)  if it has the  positive quadrant $PQD$, (negative quadrant dependence, $NQD$)dependence property.\\
\hspace*{1cm} In other words, for $n=2,$ $PUOD$ ($NUOD$) is equivalent to $PLOD$ ($NLOD$) , called as PQD (NQD) respectively.   This is given in the following definition (cf. Nelson (2006), Kotz et al. (2020)). 
\begin{d1}
\label{imp}
Let $T_1$ and $T_2$ be two random variables. $T_1$ and $T_2$ are
said to be positively (negatively) quadrant dependent denoted by $PQD$ ($NQD$) if 
$$P(T_1 \leq t_1,T_2 \leq t_2) \geq (\leq) P(T_1 \leq t_1)P(T_2 \leq t_2)$$
or equivalently 
$$P(T_1 > t_1,T_2 > t_2) \geq (\leq) P(T_1 > t_1)P(T_2 > t_2),$$ where $(t_1,t_2) \in S\subseteq R^{2}.$
\end{d1}
\hspace*{1cm}The next proposition illustrates that one can specify whether there is $OA$ or $UA$ in relative error of the two-component series system if the same is known for the relative error of the two-component parallel system and vice versa. This simplifies the computation as far as $OA$ or $UA$ in $SF$ is concerned for the parallel system, given that it is known for series system and vice versa.
\begin{p1}
Let a series and a parallel system be formed by two components. If $E_{S}^{\overline{F}}(t)$ and $E_{P}^{\overline{F}}(t)$ represent the relative error incurred in $SF$ of the series and parallel system respectively due to deviation from independence assumption then $E_{S}^{\overline{F}}(t)$ and $E_{P}^{\overline{F}}(t)$ are of opposite signs for any $t>0.$
\end{p1}
{\bf Proof.} The proof follows as a direct application of Definition \ref{imp} of PQD (LQD).$\hfill\Box$\\
\hspace*{1cm}As a result, we infer that for a given bivariate distribution, if there is an over-assessment (under-assessment) in the survival function for a two-component series system, then there is an under-assessment (over-assessment) in the survival function for a two-component parallel system. So, for a bivariate distribution one may prefer to look into whether the relative error is negative or positive for a series system and conclude that it is just the opposite for parallel systems.
\subsection{Role of stochastic orders with an important note on existing literature}
Stochastic orders are quite well known in the field of reliability theory. A brief review on stochastic orders is placed for beginners (cf. Shaked and Shanthikumar (2007)) by the following definition.
\begin{d1} Let $X$ and $Y$ be two random variables with $SF, FR, RFR, MRL,AI$ denoted by $\bar{F}_{X}, r_{X}, \mu_{X}, M_X, L_{X}$ and $\bar{F}_{Y}, r_{Y}, \mu_{Y}, m_Y, L_{Y}$ respectively. Then \begin{enumerate}
\item [$(i)$] $X$  is said to be less (greater) than $Y$ in usual stochastic order ($ST$), denoted by\\ $X\leq_{ST}(\geq_{ST})Y$ if $\bar{F}_{X}(t) \leq (\geq) \bar{F}_{Y}(t)$ for all $t\geq 0.$
\item [$(ii)$]  $X$   is said to be less (greater) than $Y$ in failure rate order ($FR$), denoted by\\ $X\leq_{FR}(\geq_{FR}~)Y$ if $r_{X}(t) \geq (\leq) r_{Y}(t)$ for all $t\geq 0.$
\item [$(iii)$]  $X$   is said to be less (greater) than $Y$ in reversed failure rate order ($RFR$), denoted by \\$X\leq_{RFR}(\geq_{RFR})Y$ if $\mu_{X}(t) \leq (\geq) \mu_{Y}(t)$ for all $t\geq 0.$
\item [$(iv)$]  $X$   is said to be less (greater) than $Y$ in mean residual life  ($MRL$), denoted by\\ $X\leq_{MRL}(\geq_{MRL})Y$ if $m_{X}(t) \leq (\geq) m_{Y}(t)$ for all $t\geq 0.$
\item [$(v)$]  $X$   is said to be less (greater) than $Y$ in likelihood ratio order  ($LR$), denoted by\\ $X\leq_{LR}(\geq_{LR})Y$ if $f_{X}(t)/f_{Y}(t)$ is non-increasing in $t$ for all $t\geq 0.$
\item [$(vi)$]  $X$   is said to be less (greater) than $Y$ in aging faster order  ($AF$), denoted by \\$X\leq_{AF}(\geq_{AF})Y$ if $r_{X}(t)/r_{Y}(t)$ is non-increasing in $t$ for all $t\geq 0.$
\item [$(vii)$]  $X$   is said to be less (greater) than $Y$ in aging intensity function order ($AI$), denoted by\\ $X\leq_{AI}(\geq_{AI})Y$ if $L_{X}(t) \geq (\leq) L_{Y}(t)$ for all $t>0.$
\end{enumerate}
\end{d1}
One can browse the literature to know more about stochastic orders. The chain of implication that exists among different stochastic orders are given in the following remark. 
\begin{r1}
\label{order}
$LR$ order implies the $ST, FR, RFR, MRL$ orders. $FR$ order implies $ST$ and $MRL$ order. $RFR$ order implies the $ST$ order. $AF$ order implies $AI$ order (cf. Nanda et al. (2007)).
\end{r1}
The over-assessment (relative error being negative) and under- assessment (relative error being positive) of different reliability measures due to the assumption of independence have been evaluated separately for different bivariate distribution by previous authors.
In view of this, we propose our findings in the following remark. To the best of our knowledge, previous authors have not taken into account the role of stochastic orderings in error analysis which arise due to deviation from dependence among system components.
\begin{r1}
\begin{enumerate}
\item[(i)] The over-assessment and under-assessment of a reliability function (say SF, FR, RFR, MRL, AI) are basically the study of the stochastic ordering (respectively, usual stochastic order, FR order, RFR order, MRL order, AI order, respectively) among the corresponding function between the dependent and independent counterpart of the system concerned. Mathematically, we express as, let $T_D$ and $T_I$ represent the lifetime of the series system formed by $n$ dependent and independent components respectively.\\
${}$\hspace{0.5cm}Clearly, $T_{D}\leq_{ST} (\geq_{ST})  T_{I}$ if and only if there is $OA$ ($UA$) in $SF$. So is for orders based on  $FR, MRL, RFR, AI,$ i.e., $T_{D}\leq_{FR} (\geq_{FR})  T_{I}$ if and only if there is $UA$ ($OA$) in $FR$., $T_{D}\leq_{MRL} (\geq_{MRL})  T_{I}$ if and only if there is $OA$ ($UA$) in $MRL$., $T_{D}\leq_{RFR} (\geq_{RFR})  T_{I}$ if and only if there is $OA$ ($UA$) in $RFR$., $T_{D}\leq_{AI} (\geq_{AI})  T_{I}$ if and only if there is $UA$ ($OA$) in $AI$. 
\item[(ii)] A stronger ordering between dependent and independent systems pertaining to a function will help us to conclude about other orderings that prevail in the stochastic ordering implication chain. As a result, one can assess the relative error of functions that follow from a stronger ordering. \\${}$\hspace{0.5cm}For example, if there is over-assessment (under-assessment) in failure rate then obviously there is under-assessment (over-assessment) in survival function, mean residual function etc.
\item[(iii)] An alternative approach could be looking at the ratio of respective densities, survivals, distributions, failure rates of systems having dependent and independent components to infer whether there is likelihood ratio order ($lr$), FR order, aging faster order, respectively.\\${}$\hspace{0.5cm}
If there is $lr$ order between systems that have dependent and independent components, then one can accordingly conclude about the assessment of other functions. 
\end{enumerate}
\end{r1}
\hspace*{1cm} The next theorem helps us to determine the sign of relative error of survival function if the sign of failure rate function is known for all $t$ or for all $t\leq t_{0}$ where $t_0$ is some positive number.
\begin{p1}
If $E_{S}^{r}(t) \leq (\geq)~ 0 $ for all $t\geq 0,$ then $E_{S}^{\overline{F}}(t)\geq (\leq) ~0$ for $t\geq 0.$$\hfill\Box$
\end{p1}
\begin{c1}
Let $t_0>0.$
If $E_{S}^{r}(t) \leq (\geq)~ 0 $ for all $t\leq t_0,$ then $E_{S}^{\overline{F}}(t)\geq (\leq)~ 0$ for $t\leq t_0.$$\hfill\Box$
\end{c1}
\section{Error Analysis of multi-component series systems}
 \hspace*{1cm}  In literature, we find the works of Gupta and Gupta (1990), Klein and Moeschberger (1986, 1987), and Moeschberger and Klein (1984) in the calculation and analysis of relative errors. Recently, Nanda et al. (2022) focused on relative errors in the various reliability measures incurred due to assumption of independence of components in a two component series system when components actually follow bivariate exponential distributions.\\ 
 \hspace*{1cm} We state a lemma from Nanda et al. (2022) to be used in upcoming section.
\begin{l1}
\label{lemma}
For $t>0,$ $$g(t)=\frac{\gamma}{\beta}\Big(\frac{e^{\beta t}-1}{e^{\gamma t}-1}\Big)-1,$$
is increasing (resp. decreasing) in $t,$ provided $\beta>\gamma$ (resp. $\beta<\gamma$).
\end{l1}
The above Lemma \ref{lemma} may be taken to arrive at the upcoming Remark \ref{impr}.
\begin{r1}
\label{impr}
Note that, for a function $g$ as defined in Lemma \ref{lemma}, $g(t)<0$ for $\gamma>\beta.$ This is because, if  $\gamma>\beta$ then $g(t)$  is decreasing in $t,$ giving $g(t)<g(0).$ Here, $g(0)=0,$ as
\begin{equation}
\lim_{t\rightarrow 0} \frac{\gamma}{\beta}\Big(\frac{e^{\beta t}-1}{e^{\gamma t}-1}\Big)=1.\nonumber
\end{equation}
\end{r1}
\subsection{Series system with components having multivariate Exponential distribution}
To begin with, we consider an $n$ component series system with dependent components that follow some multivariate exponential distribution. To compute the associated errors, we must necessarily look at a series system with $n$ independent components, each of which has an exponential distribution. Henceforth, for the sake of simplicity, in this section we denote the series system lifetime by $T$ instead of $T^{S}$ and drop $S$ (representing series system) from all the functions, viz., $SF, FR, RFR, AI.$ 
\subsubsection{Independent Exponential}
\label{indenp}
The survival function of $n$ independent random variables $T_1,T_2,\ldots,T_n$ having exponential distribution is given by \begin{equation}\label{indee}\overline{F}(t_1,t_2,\ldots,t_n)=\exp\Big\{-\sum_{i=1}^{n}\lambda_i t_i\Big\}.\end{equation}
The $SF, FR, RFR,$ and $AI$ of the series system with lifetime $T_I$ formed out of $n$ independent components each having exponential distribution with survival function given in (\ref{indee}) are obtained using (\ref{seriesi}) as
$$\overline{F}_{I}(t)=\exp\Big\{-t\sum_{i=1}^{n}\lambda_i\Big\}, r_{I}(t)=\sum_{i=1}^{n}\lambda_i, \mu_{I}(t)=\frac{\sum_{i=1}^{n}\lambda_i}{\exp\Big\{t\sum_{i=1}^{n}\lambda_i \Big\}-1},$$ for all $t\geq 0,$ and $L_{I}(t)=1,$ for all $t> 0.$
Since $r_{I}(t)=\sum_{i=1}^{n}\lambda_i,$ it is clear that the resultant series system is also exponentially distributed with parameter $\sum_{i=1}^{n}\lambda_i$ (also evident since $L_{I}(t)=1$) and has a decreasing reversed failure rate.\\
In the following discussion, we take up various forms of multivariate exponential distributions for further investigation in accordance with our present study.
\subsubsection{Marshall and Olkin's Multivariate Exponential distribution ($MOME$)}
In forthcoming discussions, we use the same notation $T_{D}$ to represent the lifetime of series system with dependent components irrespective of the multivariate distribution that it follows.  \\

The joint distribution of $T_1,T_2,\ldots,T_n$ following Marshall and Olkin's Multivariate Exponential distribution ($MOME$) is given by
\begin{eqnarray}
\overline{F}(t_1,t_2,\ldots,t_k)&=&\exp\big\{-\sum_{i=1}^{n}\lambda_i t_i-\sum_{i_1<}\sum_{i_2}\lambda_{i_1,i_2}\max(t_{i_1},t_{i_2})-\sum_{i_1<}\sum_{i_2<}\sum_{i_3}\lambda_{i_1,i_2,i_3}\max(t_{i_1},t_{i_2},t_{i_3})\nonumber\\&&
\label{momed44}\ldots-\lambda_{i_1,i_2,\ldots,i_k}\max(t_{i_1},t_{i_2},\ldots,t_{i_n})
\big\},
\end{eqnarray}
where the constants that appear in (\ref{momed44}) are all non-negative.
The $SF, FR, RFR$ and $AI$ function of a series system with lifetime, denoted by $T_D$ formed out of $n$ dependent components having $MOME$ distribution are obtained using (\ref{momed44}) and (\ref{seriesd}),
$$\overline{F}_{D}(t)=\exp(-\lambda t), r_{D}(t)=\lambda, \mu_{D}(t)=\frac{\lambda }{\exp\big(\lambda t\big)-1}, L_{D}(t)=~1,$$
where \begin{equation}
\label{lambmome}
\lambda=\Big(\sum_{i=1}^{n}\lambda_i+\sum_{i_1<}\sum_{i_2}\lambda_{i_1,i_2}+\sum_{i_1<}\sum_{i_2<}\sum_{i_3}\lambda_{i_1,i_2,i_3}+\ldots+\lambda_{i_1,i_2,\ldots,i_n}\nonumber
\Big).
\end{equation}
We note that, $$\frac{f_{D}(t)}{f_{I}(t)}=\Big(\frac{\lambda}{\sum_{i=1}^{n}\lambda_i}\Big)\exp\big(-(\lambda-\sum_{i=1}^{n}\lambda_i)t\big)$$ is decreasing in $t,$ so there is $LR$ ordering between the lifetimes of series system having dependent and independent components respectively, i.e., $T_{D}\leq_{LR} T_{I}.$ Also, $$\frac{r_{D}(t)}{r_{I}(t)}=\frac{\lambda}{\sum_{i=1}^{n}\lambda_i}, t\geq 0,$$ giving $T_{D}\stackrel{AI}{=}T_{I}.$  
This implies $T_{D}\leq_{FR}T_{I}$ (equivalent to $r_{D}(t) \geq r_{I}(t)$), $T_{D}\leq_{ST} T_{I}$ (equivalent to $\bar{F}_{D}(t) \leq \bar{F}_{I}(t)$), $T_{D}\leq_{MRL} T_{I}$ (equivalent to $m_{D}(t) \leq m_{I}(t)$) and $T_{D}\leq_{RFR} T_{I}$ (equivalent to $\mu_{D}(t) \leq \mu_{I}(t)$) for $t\geq 0,$ resulting in $OA$ in $SF, RFR, MRL$ and $UA$ in $FR.$ \\
 \hspace*{1cm} We now focus on relative errors, denoted by $E_1$ for $MOME$ distribution of individual functions (excepting $MRL,$ to avoid lengthy computation) and found that
$$E_{1}^{\overline{F}}(t)=\exp\Big\{-t\Big(\lambda-\sum_{i=1}^{n}\lambda_i\Big)\Big\}-1,$$ $$E_{1}^{r}(t)=\frac{\lambda-\sum_{i=1}^{n}\lambda_i}{\sum_{i=1}^{n}\lambda_i}, E_{1}^{\mu}(t)=\Big(\frac{\lambda}{\sum_{i=1}^{n}\lambda_i}\Big)\Big\{\frac{\exp(t\sum_{i=1}^{n}\lambda_i)-1}{\exp(\lambda t)-1}\Big\}-1, E_{1}^{L}(t)=0,$$
i.e., the relative error in $FR$  and $AI$ is a nonnegative constant. The relative error in  $SF$ and $RFR$ are decreasing in $t.$ The absolute value of the relative error is found to be the least in the $AI$ function, as expected since the $AI$ function gives an overall idea of the phenomenon of aging of a system. It is found to be least affected by the assumption of independence of components, when the components are actually dependent.  The absolute value of the relative error in the $SF$ and the $RFR$ are bounded between 0 and 1.
\subsubsection{Multivariate Gumbel's Type I ($MGI$)}
\hspace*{1cm}
We  now take up multivariate Gumbel' s Type I  exponential distribution ($MGI$).
The resulting joint survival distribution of lifetimes $T_1,T_{2}\ldots,T_n$ of the components having $MGI$ is given by
\begin{eqnarray}
\overline{F}(t_1,t_2,\ldots,t_n)&=&\exp\Big\{-\sum_{i=1}^{n}\lambda_i t_i-\sum_{i_1<}\sum_{i_2}\lambda_{i_1,i_2}t_{i_1}t_{i_2}-\sum_{i_1<}\sum_{i_2<}\sum_{i_3}\lambda_{i_1,i_2,i_3}t_{i_1}t_{i_2}t_{i_3}-\ldots-\nonumber\\
\label{mg1}
&&\lambda_{i_1,i_2,\ldots,i_n}t_{i_1}t_{i_2}\ldots t_{i_n})
\Big\}
\end{eqnarray}
where the constants plugged in (\ref{mg1}) are non-negative.\\
 \hspace*{1cm} To obtain the $SF, FR,$ and $RFR$  of a series system with lifetime $T_{D}$ formed out of $n$ dependent components having $MG1$ distribution, we use   (\ref{seriesd}) and (\ref{mg1}). Thus, 
$$\overline{F}_{D}(t)=\exp\big\{-a_1t -a_2t^2-a_3t^3-\ldots-a_nt^n
\big\},
r_{D}(t)=a_1 +2a_2t+3a_3t^2+\ldots+na_nt^{n-1}$$
\begin{equation}
\mu_{D}(t)=\frac{a_1 +2a_2t+3a_3t^2+\ldots+na_nt^{n-1}}{\exp\big\{a_1t+a_2t^2+a_3t^3+\ldots+a_nt^n
\big\}-1}=\frac{\theta^{'}(t)}{\exp(\theta(t))-1},\nonumber
\end{equation}
where $\theta(t)
=a_1t+a_{2}t^{2}+a_3t^{3}+\ldots+a_nt^{n}$ and $$a_1=\sum_{i=1}^{n}\lambda_i, a_2=\sum_{i_1<}\sum_{i_2}\lambda_{i_1,i_2}, a_3=\sum_{i_1<}\sum_{i_2<}\sum_{i_3}\lambda_{i_1,i_2,i_3};~, \ldots,~ a_n=\lambda_{i_1,i_2,\ldots,i_n}.$$ Here
$\theta^{'}(t)$ represents differentiation of $\theta(t)$ with respect to $t.$ The $AI$ function of $T_{D}$ is
\begin{equation}
\label{new1}
L_{D}(t)=\frac{a_1t+2a_2 t^2+3a_3 t^3+\ldots+na_n t^{n}
}{a_1t+a_2t^2+a_3t^3+\ldots+a_nt^n
}.
\end{equation}
\hspace*{1cm}Now, highlight a mathematical tool from Bhattacharjee et al. (2013), which is to be used in following analysis.
\begin{l1}
\label{bhatta}
Let $p_i > 0, q_i > 0,$ for $i = 1, 2,\ldots, n.$ Then
$$\min_{1\leq i\leq k}\big(\frac{p_i}{q_i}
\big) \leq \frac{\big(\sum_{1\leq i\leq n}p_i\big)}{\big(\sum_{1\leq i\leq n}q_i\big)}\leq \max_{1\leq i\leq n}\big(\frac{p_i}{q_i}
\big) $$$\hfill\Box$
\end{l1}
\hspace*{1cm}We mention a lemma from Bhattacharjee (2022) which is also established by Sunoj and Rasin (2018) using the concept of quantile function.
\begin{l1}\label{theoremss}
	A random variable $X$ is $IFRA$ (resp. $DFRA$) if and only if  $L(t) \geqslant (resp. \;\leqslant) 1$ for $t>0.$$\hfill\Box$
\end{l1}    
Applying Lemma \ref{bhatta} in (\ref{new1}), we get 
$1\leq L_{D} (t)\leq n$ and  Lemma \ref{theoremss} will imply that resultant series system is $IFRA.$ 
Here, $$\frac{\overline{F}_{D}(t)}{\overline{F}_{I}(t)}=\exp\big\{-a_2t^2-a_3t^3-\ldots-a_nt^n
\big\},$$ is decreasing in $t,$ and hence  $T_{D}\leq_{FR} T_{I}$. This in turn implies $T_{D}\leq_{ST} T_{I},$ $T_{D}\leq_{MRL} T_{I}.$ \\
Since, $$\frac{r_{D}(t)}{r_{I}(t)}=\frac{a_1 +2a_2t+3a_3t^2+\ldots+na_nt^{n-1}}{a_1},$$ is increasing in $t,$ we have $T_{D} \leq_{AFR} T_{I},$ giving $T_{D} \leq_{AI} T_{I},$ (cf. Remark \ref{order}) i.e., $L_{D}(t)\geq L_{I}(t),$ i.e., $E_{2}^{L} (t) \geq 0$ for all $t\geq 0.$\\ \hspace*{1cm} Thus, we find $OA$ in $SF, MRL$ and $UA$ in $FR, AI.$ To navigate,
we determine the relative errors in the respective functions to know their absolute value and other properties. Clearly, 
$$E_{2}^{\overline{F}}(t)=\exp\big\{ -a_2t^2-a_3t^3-\ldots-a_nt^n
\big\}-1, E_{2}^{r}(t)=\frac{2a_2t+3a_3t^2+\ldots+na_nt^{n-1}
}{a_1},$$
$$E_{2}^{\mu}(t)= \Big(\frac{\theta^{'}(t)}{a_1}\Big)\Big(\frac{e^{a_1t}-1}{e^{\theta(t)}-1}\Big)-1, E_{2}^{L}(t)=\frac{a_2 t^2+2a_3t^3+\ldots+(n-1)a_n t^{n}
}{a_1t+a_2 t^2+a_3t^3+\ldots+a_n t^n},
$$
where \begin{eqnarray}\theta(t)&=&t^2\sum_{i_1<}\sum_{i_2}\lambda_{i_1,i_2}+t^3\sum_{i_1<}\sum_{i_2<}\sum_{i_3}\lambda_{i_1,i_2,i_3}+\ldots+\lambda_{i_1,i_2,\ldots,i_n}t^n,\nonumber\\
&=&a_{1}t^{2}+a_2t^{3}+\ldots+a_nt^{n}\nonumber\end{eqnarray} and
$\theta^{'}(t)$ represents differentiation of $\theta(t)$ with respect to $t.$ \\
In $MGI$ model, relative error in $SF$ is decreasing in $t,$  for $FR$, it is increasing in $t$. The relative error in $AI$ function is decreasing in $t,$ An upper bound of relative error in $AI$ function  is $(n-1),$ (by applying Lemma \ref{bhatta}).\\
A similar study can be extended for other multivariate exponential distributions.
\subsection{Series system with components having Multivariate Weibull distribution}
In this section, we move on to multivariate Weibull distributions, which are generalizations of the multivariate exponential distribution.
We now consider a few multivariate Weibull distributions and make a similar study as discussed in earlier section.
\subsubsection{Independent Weibull}
To start with, we essentially have to figure out the reliability measures for a series system comprising of $n$ independent components having a Weibull distribution.
The $SF$ of $n$ independent random variables $T_1,T_2,\ldots,T_n$ having  Weibull distribution is given by
\begin{equation}
\label{indwe12}
\overline{F}(t_1,t_2,\ldots,t_k)=\exp\Big\{-\sum_{i=1}^{n}\lambda_i t_i^{\alpha_i}\Big\},
\end{equation}
where $\lambda_{i}>0, \alpha_{i}>0$ and $t_{i}\geq 0$ for $i=1,2,\ldots,n.$
The $SF, FR, RFR$ and $AI$,  of a series system formed out of $n$ independent components each having Weibull distribution with survival function given in (\ref{indwe12}) are obtained using (\ref{seriesi}) as
$\overline{F}_{I}(t)=\exp\Big\{-\sum_{i=1}^{k}\lambda_i t^{\alpha_i}\Big\}, r_{I}(t)=\sum_{i=1}^{k}\lambda_i \alpha_i t^{\alpha_i-1},$ and $\mu_{I}(t)=\frac{\sum_{i=1}^{n}\lambda_i \alpha_i t^{\alpha_i-1}}{\exp\Big\{\sum_{i=1}^{n}\lambda_i t^{\alpha_i}\Big\}-1}, L_{I}(t)=\frac{\sum_{i=1}^{n}\lambda_i \alpha_i t^{\alpha_i}}{\sum_{i=1}^{n}\lambda_i t^{\alpha_i}}.$ \\
\hspace*{1cm} To know the aging behavior of the resultant series system, we apply the lemma \ref{bhatta} and note that $\min_{i} \alpha_i\leq L_{I}(t)\leq \max_{i} \alpha_i$. Further using Lemma \ref{theoremss}, we claim that the system is $IFRA$ ($DFRA$) according as $\min \alpha_i>1 ~(\max\alpha_i<1).$ 
\subsubsection{Marshal and Olkin (1967a)($MOMW$)}
Marshal and Olkin (1967a) and Lee and Thompson (1974) gave a generalization of multivariate Weibull distribution, called as Marshal and Olkin multivariate Weibull Models ($MOMW$) with $SF$ given by
\begin{eqnarray}
\overline{F}(t_1,t_2,\ldots,t_n)&=&\exp\Big\{-\sum_{i=1}^{n}\lambda_i t_i^{\alpha_i}-\sum_{i_1<}\sum_{i_2}\lambda_{i_1,i_2}\max(t_{i_1}^{\alpha_{i_1}},t_{i_2}^{\alpha_{i_2}})\nonumber\\
\label{momed}
&&-\sum_{i_1<}\sum_{i_2<}\sum_{i_3}\lambda_{i_1,i_2,i_3}\max(t_{i_1}^{\alpha_{i_1}},t_{i_2}^{\alpha_{i_2}},t_{i_3}^{\alpha_{i_3}})-\nonumber\\&&\ldots-\lambda_{i_1,i_2,\ldots,i_n}\max(t_{i_1}^{\alpha_{i_1}},t_{i_2}^{\alpha_{i_2}},\ldots,t_{i_k}^{\alpha_{i_n}})
\Big\}\nonumber
\end{eqnarray}
The $SF, FR, RFR$ and $AI$ of a series system formed out of $n$ components having $MOMW$ are respectively given by
$\overline{F}_{D}(t)=\exp(A(t)),$ $r_{D}(t)=-A^{'}(t),$
$\mu_{D}(x)=\frac{A^{'}(x)}{1-\exp(-A(x))},$ $L_{D}(x)=\frac{xA^{'}(x)}{A(x)}$
where 
\begin{eqnarray}
  A(t)&=&-\Big\{\sum_{i=1}^{n}\lambda_i t^{\alpha_i}+\sum_{i_1<}\sum_{i_2}\lambda_{i_1,i_2}~t^{\max(\alpha_{i_1}, \alpha_{i_2})}+\ldots+\lambda_{i_1,i_2,\ldots,i_n}t^{\max(\alpha_{i_1}, \alpha_{i_2},\ldots,\alpha_{i_n})}
\Big\}  \nonumber\\
&=&-\Big\{\sum_{i=1}^{n}\lambda_i t^{\alpha_i}+B(t)\}\nonumber
\end{eqnarray}
where $B(t)=\sum_{i_1<}\sum_{i_2}\lambda_{i_1,i_2}~t^{\max(\alpha_{i_1}, \alpha_{i_2})}+\ldots+\lambda_{i_1,i_2,\ldots,i_n}t^{\max(\alpha_{i_1}, \alpha_{i_2},\ldots,\alpha_{i_n})}$ (say), respectively.\\
Clearly $$\frac{\overline{F}_{D}(t)}{\overline{F}_{I}(t)}=\exp\big(-B(t)\big)$$ is decreasing in $t \geq 0$, and hence $r_{D}(t)\geq r_{I}(t)$ giving $\overline{F}_{D}(t)\leq \overline{F}_{I}(t)$ and $m_{D}(t)\leq m_{I}(t)$ for all $t\geq 0.$ Thus, the relative errors in $SF, FR, MRL$ respectively satisfy $E_{3}^{\overline{F}}(t) \leq 0,$ $E_{3}^{r}(t) \geq 0,$ and $E_{3}^{m}(t) \leq 0$ for all $t\geq 0.$ In other words, there is $OA$ in $SF, MRL$ and $UA$ in $FR.$\\

 The relative error incurred in $SF, FR, RFR$ and $AI$ functions by assuming independence of components where the components are actually dependent having $MOMW$ distribution in formation of a series system are respectively given by
$$E_{3}^{\overline{F}}(t)=\exp\Big(A(t)+\sum_{i=1}^{n}\lambda_i t^{\alpha_i}\Big)-1=\exp(-B(t))-1,$$
$$E_{3}^{r}(t)=\left(-A^{'}(t)-\sum_{i=1}^{n}\lambda_i\alpha_{i} t^{\alpha_i-1}\right)\Big\{\sum_{i=1}^{n}\lambda_i \alpha_i t^{\alpha_i-1}\Big\}^{-1},$$
\begin{eqnarray}
E_{3}^{\mu}(t)=\frac{\mu_{D}(t)}{\mu_{I}(t)}-1&=&  \frac{-\Big\{\exp(\sum_{i=1}\lambda_{i}t^{\alpha_{i}})-1\Big\}\Big\{\sum_{i=1}^{n}\lambda_i \alpha_i t^{\alpha_i-1}+B^{'}(t)\Big\}}{\Big\{\sum_{i=1}^{n}\lambda_i \alpha_i t^{\alpha_i-1}\Big\}\Big\{1-\exp(\sum_{i=1}\lambda_{i}t^{\alpha_{i}})\exp(B(t))\Big\}}-1 \nonumber\\
&=&\Big(1+\frac{B^{'}(t)}{\sum_{i=1}^{n}\lambda_i \alpha_i t^{\alpha_i-1}}\Big)\Big\{\frac{1-\exp(\sum_{i=1}\lambda_{i}t^{\alpha_{i}})}{1-\exp(\sum_{i=1}\lambda_{i}t^{\alpha_{i}})\exp(B(t))}\Big\}-1,\nonumber
\end{eqnarray}
 \begin{eqnarray}
E_{3}^{L}(t)=\frac{L_{D}(t)}{L_{I}(t)}-1&=&  \frac{t\Big\{\sum_{i=1}\lambda_{i}t^{\alpha_{i}}\Big\}\Big\{\sum_{i=1}^{n}\lambda_i \alpha_i t^{\alpha_i-1}+B^{'}(t)\Big\}}{\Big\{\sum_{i=1}^{n}\lambda_i \alpha_{i}t^{\alpha_i}\Big\}\Big\{\sum_{i=1}\lambda_{i}t^{\alpha_{i}}+B(t)\Big\}}-1 \nonumber\\
&=&\Big(1+\frac{tB^{'}(t)}{\sum_{i=1}^{n}\lambda_i \alpha_i t^{\alpha_i}}\Big)\Big\{\frac{\sum_{i=1}\lambda_{i}t^{\alpha_{i}}}{\sum_{i=1}\lambda_{i}t^{\alpha_{i}}+B(t)}\Big\}-1,\nonumber
\end{eqnarray}
and 
$$E_{3}^{L}(t)=\left(-A^{'}(t)-\sum_{i=1}^{n}\lambda_i\alpha_{i} t^{\alpha_i-1}\right)\Big\{\sum_{i=1}^{n}\lambda_i \alpha_i t^{\alpha_i-1}\Big\}^{-1}.$$
$E_{3}^{L}(t)$ changes sign.
\subsubsection{Multivariate Crowder Weibull (1989) ($MCW$)}
\hspace*{1cm}Hougaard (1986,1989) presented a multivariate  Weibull distribution with joint survival function
\begin{equation}
\overline{F}(t_1,t_2,\ldots,t_n)=\exp\Big\{-\big(\sum_{i=1}^{n}\lambda_i t_{i}^{\alpha_i}\big)^{l}\Big\},\nonumber
\end{equation}
where $l>0, \lambda_i,\alpha_i \geq 0$ and $t_i \geq 0.$ Crowder (1989) extended Hougaard's distributions and proposed multivariate distributions with Weibull connection with
\begin{equation}
\label{crow}
\overline{F}(t_1,t_2,\ldots,t_n)=\exp\Big\{\gamma^{l}-\big(\gamma+\sum_{i=1}^{n}\lambda_i t_{i}^{\alpha_i}\big)^{l}\Big\}
\end{equation}
where $l>0,  \gamma \geq 0$ and $\lambda_i,\alpha_i >0.$ If $\gamma=0, l=1$ then (\ref{crow}) reduces to independent Weibull distribution. In the special case where $\alpha_1=\ldots=\alpha_n=\alpha,$ the marginals are all Weibull with the same parameter. \\
\hspace*{1cm}The $SF, FR, RFR$ and $AI$ of a series system formed out of $n$ components having $MCW$ are respectively given by
$$\overline{F}_{D}(t)=\exp\Big\{\gamma^{l}-\big(\gamma+\sum_{i=1}^{n}\lambda_i t^{\alpha_i}\big)^{l}\Big\}, r_D(t)= l\Big(\gamma+\sum_{i=1}^{n}\lambda_i t^{\alpha_i}\Big)^{l-1}\Big(\sum_{i=1}^{n}\lambda_i \alpha_i t^{\alpha_i-1}\Big),$$
$$\mu_D(t)=\frac{l\Big(\gamma+\sum_{i=1}^{n}\lambda_i t^{\alpha_i}\Big)^{l-1}\Big(\sum_{i=1}^{k}\lambda_i \alpha_i t^{\alpha_i-1}\Big)}{\exp\Big\{-\gamma^{l}+\big(\gamma+\sum_{i=1}^{n}\lambda_i t^{\alpha_i}\big)^{l}\Big\}-1},$$ for $t\geq 0$ and
$$L_D(t)=\frac{l\Big(\gamma+\sum_{i=1}^{n}\lambda_i t^{\alpha_i}\Big)^{l-1}\Big(\sum_{i=1}^{n}\lambda_i \alpha_i t^{\alpha_i}\Big)}{\exp\Big\{\gamma^{l}-\big(\gamma+\sum_{i=1}^{n}\lambda_i t^{\alpha_i}\big)^{l}\Big\} \Big\{\gamma^{l}-\big(\gamma+\sum_{i=1}^{n}\lambda_i t^{\alpha_i}\big)^{l}\Big\}}, t>0.$$
  The relative error incurred in $SF,FR, RFR$ and $AI$ by assuming independence of components where the components are actually dependent in formation of a series system are respectively given by
 $$E_{4}^{\overline{F}}(t)=\exp\big\{\gamma^{l}+\sum_{i=1}^{n} \lambda_{i}t^{\alpha_i}-\big(\gamma+\sum_{i=1}^{n} \lambda_{i}t^{\alpha_i}\big)^{l}\big\}-1,E_{4}^{r}t)=l\Big(\gamma+\sum_{i=1}^{n} \lambda_{i}t^{\alpha_i}\Big)^{l-1}-1,$$
 $$E_{4}^{\mu}(t)=\frac{l\Big(\gamma+\sum_{i=1}^{n}\lambda_i t^{\alpha_i}\Big)^{l-1}}{\Big\{\exp\big\{-\gamma^{l}+\big(\gamma+\sum_{i=1}^{n}\lambda_i t^{\alpha_i}\big)^{l}\big\}-1\Big\}}\Big\{\exp\Big(\sum_{i=1}^{n}\lambda_i t^{\alpha_i}\Big)-1\Big\}.$$
 $$E_{4}^{L}(t)=\frac{l\Big(\gamma+\sum_{i=1}^{n}\lambda_i t^{\alpha_i}\Big)^{l-1}\Big(\sum_{i=1}^{n}\lambda_i t^{\alpha_i}\Big)}{\Big\{\exp\Big(\gamma^{l}-\big(\gamma+\sum_{i=1}^{n}\lambda_i t^{\alpha_i}\big)^{l}\Big)\Big\} \Big\{\gamma^{l}-\big(\gamma+\sum_{i=1}^{n}\lambda_i t^{\alpha_i}\big)^{l}\Big\}}-1.$$
In particular, Lee (1979): II ($MLII$)
\begin{equation}
\overline{F}(x_1,x_2,\ldots,x_k)=\exp\big\{-\big(\sum_{i=1}^{k} \lambda_{i}t_i^{\alpha_i}\big)^{l}\big\}\nonumber
\end{equation}
where $\alpha_i>0, 0<\gamma\leq 1, \lambda_i>0,t_i\geq 0.$
\subsubsection{Multivariate Lee (1979) ($ML$)}
\begin{eqnarray}
\overline{F}(t_1,t_2,\ldots,t_n)&=&\exp\Big\{-\sum_{i=1}^{n}\lambda_i c_i^{\alpha} t_i^{\alpha}-\sum_{i_1<}\sum_{i_2}\lambda_{i_1,i_2}\max(c_{i_1}^{\alpha}t_{i_1}^{\alpha},c_{i_2}^{\alpha}t_{i_2}^{\alpha})\nonumber\\
\label{momed}
&&-\sum_{i_1<}\sum_{i_2<}\sum_{i_3}\lambda_{i_1,i_2,i_3}\max(c_{i_1}^{\alpha}t_{i_1}^{\alpha},c_{i_2}^{\alpha}t_{i_2}^{\alpha},c_{i_3}^{\alpha}t_{i_3}^{\alpha})-\nonumber\\
&&\ldots-\lambda_{i_1,i_2,\ldots,i_k}\max(c_{i_1}^{\alpha}t_{i_1}^{\alpha},c_{i_2}^{\alpha}t_{i_2}^{\alpha},\ldots,c_{i_n}^{\alpha}t_{i_n}^{\alpha})
\Big\}\nonumber
\end{eqnarray}
Taking its independent counterpart as $\overline{F}(t_1,t_2,\ldots,t_n)=\exp\Big\{-\sum_{i=1}^{n}\lambda_i c_i^{\alpha} t_i^{\alpha}\Big\},$ we note that the survival function, failure rate function, reversed failure rate and ageing intensity function of a series system with $k$ independent components are given by $\overline{F}_{I}(t)=\exp\Big\{-t^{\alpha}\sum_{i=1}^{n}\lambda_i c_i^{\alpha} \Big\},$ $r_{I}(t)=\alpha t^{\alpha-1}\sum_{i=1}^{n}\lambda_i c_{i}^{\alpha}, \mu_{I}(t)=\frac{\alpha t^{\alpha-1}\sum_{i=1}^{n}\lambda_i c_{i}^{\alpha}}{\exp\Big\{t^{\alpha}\sum_{i=1}^{n}\lambda_i c_{i}^{\alpha} \Big\}-1},$ and $L_{I}(t)=\alpha.$\\
\hspace*{1cm}The survival function, failure rate, reversed failure rate and aging intensity functions of a series system formed out of $n$ components having multivariate $ML$ distribution are respectively given by
$$\overline{F}(t)=\exp\Big\{-\lambda_{L}t^{\alpha}
\Big\}, r(t)=\alpha t^{\alpha-1}\lambda_{L},
\mu(t)=\frac{\alpha t^{\alpha-1}\lambda_{L}}{e^{\lambda_{L}t^{\alpha}}-1},
  L(t)=\alpha,$$ where \begin{eqnarray}
\lambda_{L}&=&\Big(\sum_{i=1}^{n}\lambda_i c_i^{\alpha}+\sum_{i_1<}\sum_{i_2}\lambda_{i_1,i_2}\max(c_{i_1}^{\alpha},c_{i_2}^{\alpha})+\sum_{i_1<}\sum_{i_2<}\sum_{i_3}\lambda_{i_1,i_2,i_3}\max(c_{i_1}^{\alpha},c_{i_2}^{\alpha},c_{i_3}^{\alpha})+\nonumber\\
&&\ldots+\lambda_{i_1,i_2,\ldots,i_n}\max(c_{i_1}^{\alpha},c_{i_2}^{\alpha},\ldots,c_{i_n}^{\alpha})\Big).\nonumber
\end{eqnarray}
Since $\frac{f_{D}(t)}{f_{I}(t)}$ is decreasing in $t,$ we have $T_{D}\leq_{lr} T_{I}$ and hence $T_{D}\leq_{hr} T_{I}$, $T_{D}\leq_{st} T_{I}$, $T_{D}\leq_{rhr} T_{I}$ and $T_{D}\leq_{mrl} T_{I}.$  We note that $\frac{r_{D}(t)}{r_{I}(t)}$ is non-negative constant.
The relative error in  each of survival function, reversed failure rate,  mean residual life function  is negative whereas the relative error in failure rate is positive, for aging intensity function, it is zero.\\
We give a modification of Lemma  4.2 of et al. Nanda (2022) without proof.
\begin{l1}
\label{lemma2}
For $t>0,$ $\beta,\gamma,\alpha>0,$ $$h(t)=\frac{\gamma}{\beta}\Big(\frac{e^{\beta t^{\alpha}}-1}{e^{\gamma t^{\alpha}}-1}\Big)-1,$$
is increasing (resp. decreasing) in $t,$ provided $\beta>\gamma$ (resp. $\beta<\gamma$).
\end{l1}

   {\bf Proof.} Taking $x^{\alpha}=t (say),$ and noting that $$\frac{d}{dx}h(x)=\frac{dh}{dt}\frac{dt}{dx},$$ we complete the proof.
\begin{l1}
\label{pro}
Note that, for $h$ as defined in Lemma \ref{lemma2}, $h(x)<0$ for $\gamma>\beta.$ This is because, if  $\gamma>\beta$ then $h(x)$  is decreasing in $x,$ giving $h(x)<h(0).$ Here, $h(0)=0,$ as
\begin{equation}
\lim_{x\rightarrow 0} \frac{\gamma}{\beta}\Big(\frac{e^{\beta x^{\alpha}}-1}{e^{\gamma x^{\alpha}}-1}\Big)=1\nonumber
\end{equation}
\end{l1}

 The relative error incurred in survival function, failure rate, reversed failure rate and aging intensity functions by assuming independence of components where the components are actually dependent in formation of a series system are respectively given by
$$E_{6}^{\overline{F}}(t)=\exp\Big\{-t^{\alpha}\Big(\lambda_{L}-\sum_{i=1}^{n}\lambda_{i}c_{i}^{\alpha}\Big)\Big\}-1,
E_{6}^{r}(t)=\frac{\lambda_{L}}{\sum_{i=1}^{n}\lambda_{i}c_{i}^{\alpha}}-1,$$ $$E_{6}^{\mu}(t)=\Big(\frac{\lambda_{L}}{\sum_{i=1}^{n}\lambda_ic_{i}^{\alpha}}\Big)\Big\{\frac{\exp(t^{\alpha}\sum_{i=1}^{n}\lambda_i c_{i}^{\alpha})-1}{\exp(\lambda_{L} t^{\alpha})-1}\Big\}-1,E^{L}_{6}(t)=0.$$
We conclude that the relative errors in all aforementioned functions are decreasing in $t$ except that for failure rate and AI function it is non-negative constant.
\subsubsection{Lu and Bhattcharyya (1990): I}

    Lu and Bhattacharyya I\\
As an application, Hougaard(1986) used the above Weibull model on a data on tumour appearance in 50 litters of female rats. Each litter contained one drug treated and two control rats. The data was studied to find estimate of failures in marginal distribution, maximized likelihood function and standard errors. The true maximum likelihood estimate of log relative risk was found to be 0.944 for the dependence model and 0.04 less for the independence model. The former model had standard error 0.327 whereas the later had 0.01 less. The likelihood ratio test statistic was estimated to be 7.99 for the former and 0.29 more for the later model. It was also concluded that the two stage estimate and the maximum likelihood estimate differed marginally.\\

\begin{equation}
\overline{F}(t_1,t_2,\ldots,t_k)=\exp\Big\{-\Big(\sum_{i=1}^{n}\lambda_{i}t_i^{\alpha_i}+\delta w(t_1,t_2,\ldots,t_n)\Big)\Big\}\nonumber
\end{equation}
where $w(t_1,t_2,\ldots,t_n)=\Big\{\sum_{i=1}^{n}\lambda_{i}^{\frac{1}{m}}t_i^{\frac{\alpha_i}{m}}\Big\}^{m}.$
The survival function, failure rate, reversed failure rate and aging intensity functions of a series system formed out of $n$ components having Crowder distribution are respectively given by
\begin{equation}
\overline{F}^{T}_{1}(t)=\exp\Big\{-\Big\{\sum_{i=1}^{n}\lambda_{i}t^{\alpha_i}+\delta \Big(\sum_{i=1}^{n}\lambda_{i}^{\frac{1}{m}}t^{\frac{\alpha_i}{m}}\Big)^{m}\Big\}\Big\}=\exp(A(t)-\delta B(t))\nonumber
\end{equation}
\begin{equation}
r^{T}_1(t)= \sum_{i=1}^{n}\lambda_{i}\alpha_i t_i^{\alpha_i}+\frac{\delta}{x}\Big(\sum_{i=1}^{n}\lambda_{i}^{\frac{1}{m}}\alpha_i t^{\frac{\alpha_i}{m}}\Big)\Big(\sum_{i=1}^{n}\lambda_{i}^{\frac{1}{m}}t^{\frac{\alpha_i}{m}}\Big)^{m}=-C^{'}(t)\nonumber
\end{equation}
where $C(t)=A(t)-\delta B(t),$ and $B(t)=\Big(\sum_{i=1}^{n}\lambda_{i}^{1/m}t^{\alpha_{i}/m}\Big)^{m}.$
 $$\mu^{T}_1(t)=\frac{\gamma \big(\sum_{i=1}^{n} \lambda_{i}t^{\alpha_i}\big)^{\gamma-1} \Big(\sum_{i=1}^{n}\lambda_i \alpha_i t^{\alpha_i-1}\Big)}{\exp\Big\{\sum_{i=1}^{n}\lambda_i t^{\alpha_i}\Big\}-1}=\frac{C^{'}(t)}{1-\exp(-C(t))}$$
$$L^{T}_1(t)=\frac{\Big(\sum_{i=1}^{n}\lambda_{i}\alpha_i t_i^{\alpha_i}\Big)+\delta \Big(\sum_{i=1}^{n}\lambda_{i}^{\frac{1}{m}}\alpha_i t^{\frac{\alpha_i}{m}}\Big)\Big(\sum_{i=1}^{n}\lambda_{i}^{\frac{1}{m}} t^{\frac{\alpha_i}{m}}\Big)^{m-1}}{\Big(\sum_{i=1}^{n}\lambda_{i} t_i^{\alpha_i}\Big)+\delta \Big(\sum_{i=1}^{n}\lambda_{i}^{\frac{1}{m}} x^{\frac{\alpha_i}{m}}\Big)^{m}}=-\frac{tC^{'}(t)}{C(t)}$$
$$E^{\overline{F}}(t)=\delta \exp(B(t))-1, E^{r}(t)=-\frac{\delta B^{'}(t)}{A^{'}(t)}, E^{L}(t)=\frac{\delta \Big(\frac{B^{'}(t)}{A^{'}(t)}+\frac{B(t)}{A(t)}-2\Big)}{1-\delta \frac{B(t)}{A(t)} }$$
\subsubsection{Farlie-Gumbel-Morgenstern Weibull (FGMW) System of distribution}
A n-dimensional random variable $(T_{1},T_{2},\ldots,T_{n})$ is said to follow Farlie-Gumbel-Morgenstern Weibull (FGMW) distribution if its survival function is given by
\begin{equation}
\overline{F}(t_1,t_2,\ldots,t_n)=\Big\{\prod_{i=1}^{n}\overline{F}_{i}(t_i)\Big\}\Big\{1+\gamma\prod_{i=1}^{n} [1-\overline{F}_{i}(t_i)]\Big\}\nonumber
\end{equation}
where $\overline{F}_{i}(t_i)=\exp(-\lambda_i t_{i}^{\alpha_i})$ for $i=1,2,\ldots,n$ is the survival function of $T_{i}$ and $-1<\gamma<1.$ \\
The survival function of a series system formed out of $n$ components having $FGMW$ are respectively given by
\begin{eqnarray}
\overline{F}^{T}_{1}(t)&=&\Big(\prod_{i=1}^{k}\exp(-\lambda_{i}x^{\alpha_{i}})\Big)\Bigg(1+\gamma\Big(\prod_{i=1}^{k}\big(1-\exp(-\lambda_{i}t^{\alpha_{i}})\big)\Big)\Bigg)\nonumber\\
&=&\exp(A(t))\Big(1+\gamma (1-\theta(t))\Big)=\exp(A(t))\phi(t),\nonumber
\end{eqnarray}
where $$\theta(t)=\sum_{k=1}^{n}(-1)^{k-1}\sum_{i_{1}<} \sum_{i_{2}<} \ldots\sum_{<i_{k}}\exp\Big(\sum_{p=1}^{k}(-\lambda_{i_{p}}t^{\alpha_{i_{p}}})\Big),$$ $\phi(t)=1+\gamma (1-\theta(t)),$ and $A(t)=\sum_{i=1}^{k}(-\lambda_{i}t^{\alpha_{i}}).$ It is worthwhile to immediately note that $ (1-\theta(t))$  is non-decreasing in $t$ and lies between 0 and 1. Interestingly, $\phi(t)\geq 0.$ As a result, $\phi(t)=1+\gamma (1-\theta(t))$ is increasing (decreasing) in $t$ if $\gamma \geq (\leq) 0.$  Since, $r_{I}(t)=-A^{'}(t),$ it follows that $A^{'}(t)\leq 0$ for all $t>0.$
The corresponding failure rate, reversed failure rate and aging intensity functions are given by
$r^{T}_1(x)=-A^{'}(t)-\frac{\phi^{'}(t)}{\phi(t)},$
$$\mu^{T}_1(t)=\frac{-A^{'}(t)\phi(t)-\phi^{'}(t)}{1-\phi(t)\exp(A(t))}, L^{T}_1(t)=t\Big\{\frac{ A^{'}(t)\phi(t)+\phi^{'}(t)}{\phi(t)\Big(A(t)+\ln \phi(t)\Big)}\Big\},$$
Here, $E^{\overline{F}}(t)=\gamma(1-\theta(t)), E^{r}(t)=\frac{\phi^{'}(t)}{\phi(t)A^{'}(t)},$ So,  $E^{\overline{F}}(t)$ takes same sign as $\gamma.$ Further, $E^{\overline{F}}(t)\leq (\geq) \gamma$ according as $\gamma \geq (\leq) 0.$ $E^{r}(t)$ takes opposite sign to $\gamma.$

\subsubsection{Lu and Bhattacharyya (1990): II}
\begin{equation}
\overline{F}(x_1,x_2)=\Big[1+\big[\big\{\exp[\lambda_{1}x_1^{\alpha_1}]-1\big\}^{1/\gamma}+   \big\{\exp[\lambda_{2}x_2^{\alpha_2}]-1\big\}^{1/\gamma}\big]^{\gamma}\Big]^{-1}
\end{equation}
for $-1<\gamma<1.$

\begin{equation}
\overline{F}(x_1,x_2,\ldots,x_{i})=\Big[1+\sum_{i=1}^{n}\big[\big\{\exp[\lambda_{i}t_i^{\alpha_i}]-1\big\}^{1/\gamma}\big]^{\gamma}\Big]^{-1}
\end{equation}
for $-1<\gamma<1.$

  \section{Reliability Analysis of $n$ component parallel system and Incurred Error}
\subsection{Survival function of parallel system in terms of joint survival function of lifetime multivariate distributions}
The joint distribution of any random vector $(X_1,X_2,\ldots,X_n)$ can be expressed as
\begin{equation*}
	\begin{array}{rll}
		F(x_1,x_2,\cdots,x_n)=&P(X_1\leq x_1,X_2\leq x_2,\cdots, X_n\leq x_n)=P(\cap_{i=1}^n t_i\leq t_i)&\\
		=&1-P(\cup_{i=1}^n t_i>t_i)\text{~~~(By de Morgan's Law)}&\\
		=&1-\sum_{i=1}^n \bar{F}_{t_i}(t_i)+\underset{1\leq i<j\leq n}{\sum\sum}\bar{F}_{t_i,X_j}(t_i,x_j)-\underset{1\leq i<j<k\leq n}{\sum\sum\sum}\bar{F}_{t_i,X_j,X_k}(t_i,x_j,x_k)\\
		&\ldots+(-1)^n\bar{F}(x_1,x_2,x_3,\cdots,x_n) \text{~~~(By Poincare's Theorem)}&\\
		=&1-\sum_{i=1}^n \bar{F}(-\infty,\cdots,-\infty,\underbrace{t_i}_{i^\text{th}\text{ position}},-\infty,\cdots,-\infty)&\\
		&\multicolumn{2}{l}{+\underset{1\leq i<j\leq n}{\sum\sum}\bar{F}(-\infty,\cdots,-\infty,\underbrace{t_i}_{i^\text{th}\text{ position}},-\infty,\cdots,-\infty,\underbrace{x_j}_{j^\text{th}\text{ position}},-\infty,\cdots,-\infty)}\\
		&\multicolumn{2}{l}{-\underset{1\leq i<j<k\leq n}{\sum\sum\sum}\bar{F}(-\infty,\cdots,-\infty,\underbrace{t_i}_{i^\text{th}},-\infty,\cdots, -\infty,\underbrace{x_j}_{j^\text{th}},-\infty,\cdots,-\infty,\underbrace{x_k}_{k^\text{th}},-\infty,\cdots,-\infty)}\\
		&+\ldots+(-1)^n\bar{F}(x_1,x_2,x_3,\cdots,x_n)&\\
	\end{array}
\end{equation*}
Hence, the survival function of lifetime $T$ of a parallel system is
\begin{equation*}
	\begin{array}{rll}
		 \overline{F}_{T}(t)=&1-F(t,t,t,\cdots,t)=\sum_{i=1}^n \bar{F}(0,\cdots,0,\underbrace{t}_{i^\text{th}\text{ position}},0,\cdots,0)&\\
		&\multicolumn{2}{l}{-\underset{1\leq i<j\leq n}{\sum\sum} \overline{F}(0,\cdots,0,\underbrace{t}_{i^\text{th}\text{ position}},0,\cdots,0,\underbrace{t}_{j^\text{th}\text{ position}},0,\cdots,0)}\\
		&\multicolumn{2}{l}{+\underset{1\leq i<j<k\leq n}{\sum\sum\sum}\bar{F}(0,\cdots,0,\underbrace{t}_{i^\text{th}\text{ position}},0,\cdots, 0,\underbrace{t}_{j^\text{th}\text{ position}},0,\cdots,0,\underbrace{t}_{k^\text{th}\text{ position}},0,\cdots,0)}\\
		&\vdots&\\
		&+(-1)^{n-1} \overline{F}(t,t,t,\cdots,t)&\\
	\end{array}
\end{equation*}
Note that the lower limits of $-\infty$ have been replaced by $0$ because we are dealing with lifetime variables.
In this study, we first consider an $k$ component parallel system with dependent components having Marshall and Olkin's multivariate exponential distribution.

\subsection{Parallel system with components having Multivariate Exponential distribution}
\subsubsection{Independent exponential}
\begin{equation}\label{inde}\overline{F}(x_1,x_2,\ldots,x_n)=\exp\Big\{-\sum_{i=1}^{n}\lambda_i t_i\Big\}.\nonumber
\end{equation}
\begin{equation}\label{Pindexp}
	\begin{split}
		\bar{F}_{T}(t)=&\sum_{i=1}^n\exp\left\{-t\lambda_i\right\}-
		\mathop{\sum\sum}_{1\leq i<j\leq n}\exp\left\{-t(\lambda_i+\lambda_j)\right\}+\mathop{\sum\sum\sum}_{1\leq i<j<k\leq n}\exp\left\{-t(\lambda_i+\lambda_j+\lambda_k)\right\}\\
		&+\ldots+(-1)^{k-1}\mathop{\sum\cdots\sum}_{1\leq i_1<\cdots<i_k\leq n}\exp\left\{-t\left(\sum_{p=1}^k\lambda_{i_p}\right)\right\}+\ldots+(-1)^{n-1}\exp\left\{-t\sum_{i=1}^{n}\lambda_i\right\}\\
		=&\sum_{k=1}^n(-1)^{k-1}\mathop{\sum\cdots\sum}_{1\leq i_1<\cdots<i_k\leq n}\exp\left\{-t\left(\sum_{p=1}^k\lambda_{i_p}\right)\right\}\nonumber
	\end{split}
\end{equation}
This survival can alternatively be calculated as $$\bar{F}_T(t)=1-P[\max\{T_1,T_2,\ldots,T_n\}\leq t]=1-\prod_{i=1}^n P[T_i\leq t]=1-\prod_{i=1}^n(1-e^{-\lambda_i t}),$$ but the expression obtained in \ref{Pindexp} will be rather helpful while finding the relative error.

Therefore, $$f_T(t)=-\frac{\mathrm{d}}{\mathrm{dt}}\bar{F}_T(t)=\sum_{k=1}^n(-1)^{k-1}\mathop{\sum\cdots\sum}_{1\leq i_1<\cdots<i_k\leq n}\left[\left(\sum_{p=1}^k\lambda_{i_p}\right)\exp\left\{-t\left(\sum_{p=1}^k\lambda_{i_p}\right)\right\}\right]$$
\subsubsection{Marshall and Olkin's Multivariate Exponential distribution ($MOME$)}
Kotz et al. (2000) described $MOME$ in the following words: In a system of $k$ components, the distribution of lifetimes between "fatal shocks" to the combination $\{a_1,\ldots,a_l\}$ of components is supposed to have an exponential distribution with mean $1/\lambda_{a_1,\ldots,a_l}.$ The $2^{k-1}-1$ different distributions of this kind are supposed to be a mutually independent set.
The resulting joint distribution of lifetimes $X_1,\ldots,X_n$ of the components having $MOME$ is given by
\begin{eqnarray}
\overline{F}(x_1,x_2,\ldots,x_n)&=&\exp\Big\{-\sum_{i=1}^{k}\lambda_i t_i-\sum_{i_1<}\sum_{i_2}\lambda_{i_1,i_2}\max(x_{i_1},x_{i_2})\nonumber\\&&-\sum_{i_1<}\sum_{i_2<}\sum_{i_3}\lambda_{i_1,i_2,i_3}\max(x_{i_1},x_{i_2},x_{i_3})-\nonumber\\
\label{momed4}&&\ldots-\lambda_{1,2,\ldots,n}\max(x_{1},x_{2},\ldots,x_{n})
\Big\}\nonumber
\end{eqnarray}
\begin{equation*}
	\begin{split}
		\bar{F}_{T}(t)=&\sum_{i=1}^n\exp\left\{-t\lambda_i-\sum_{j\neq i}t\lambda_{i,j}-\mathop{\mathop{\sum\sum}_{j\neq i\neq k}}_{j<k}t\lambda_{i,j,k}-\cdots-t\lambda_{1,\cdots,k}\right\}\\
		-&\mathop{\sum\sum}_{1\leq i_1<i_2\leq n}\exp\left\{-t\lambda_{i_1}-t\lambda_{i_2}-t\lambda_{i_1,i_2}-\sum_{j\notin\{i_1,i_2\}}\left(t\lambda_{i_1,j}+t\lambda_{i_2,j}\right)-\sum_{j\notin\{i_1,i_2\}}t\lambda_{i_1,i_2,j}\right.\\
		&\left. -\mathop{\sum\sum}_{j<k\notin\{i_1,i_2\}}\left(t\lambda_{i_1,j,k}+t\lambda_{i_2,j,k}\right)-\cdots-t\lambda_{1,2,\cdots,n}\right\}\\
		\vdots&\\
		+&(-1)^{n-1}\exp\left\{-t\left(\sum_{i=1}^n\lambda_i +\mathop{\sum\sum}_{i<j}\lambda_{i,j}+\mathop{\sum\sum\sum}_{i<j<k}\lambda_{i,j,k}+\cdots+\lambda_{1,2,\cdots,n}\right)\right\}\nonumber
	\end{split}
\end{equation*}
\begin{equation*}
	\begin{split}
		\bar{F}_T(t)=&\sum_{i=1}^n\exp\left\{-t\left(\lambda_i+\sum_{j\neq i}\lambda_{i,j}+\mathop{\mathop{\sum\sum}_{j\neq i\neq k}}_{j<k}\lambda_{i,j,k}+\cdots+\lambda_{1,2,\cdots,n}\right)\right\}\\
		-&\mathop{\sum\sum}_{1\leq i_1<i_2\leq n}\exp\left\{-t\left(\lambda_{i_1}+\lambda_{i_2}+\lambda_{i_1,i_2}+\sum_{j\notin\{i_1,i_2\}}(\lambda_{i_1,j}+\lambda_{i_2,j}+\lambda_{i_1,i_2,j})\phantom{\mathop{\mathop{\sum\sum}_{j,k\notin\{i_1,i_2\}}}_{j<k}}\right.\right.\\		&\left.\left.+\mathop{\mathop{\sum\sum}_{j,k\notin\{i_1,i_2\}}}_{j<k}(\lambda_{i_1,j,k}+\lambda_{i_2,j,k}+\lambda_{i_1,i_2,j,k})+\cdots+\lambda_{1,2,\cdots,n}\right)\right\}+\\
		&\ldots
		+\ldots+(-1)^{n-1}\exp\left\{-t\left(\sum_{i=1}^n\lambda_i +\mathop{\sum\sum}_{i<j}\lambda_{i,j}+\mathop{\sum\sum\sum}_{i<j<k}\lambda_{i,j,k}+\cdots+\lambda_{1,2,\cdots,n}\right)\right\}\\
	\end{split}
\end{equation*}
$$=\sum_{k=1}^n(-1)^{k-1}\mathop{\sum\cdots\sum}_{1\leq i_1<\cdots<i_k\leq n}\exp\left\{-t\left(\sum_{p=0}^k\left(\mathop{\mathop{\sum\cdots\sum}_{1\leq j_1<\cdots<j_p\leq n}}_{j_1,\cdots,j_p\notin\{i_1,\cdots,i_k\}}\sum_{q=1}^k\left(\mathop{\sum\cdots\sum}_{1\leq u_1<\cdots<u_q\leq k}\lambda_{i_{u_1},\cdots,i_{u_q},j_1,\cdots,j_p}\right)\right)\right)\right\}$$
Hence, defining $$A_{i_1,i_2,\ldots,i_k}=\sum_{p=0}^k\left(\mathop{\mathop{\sum\cdots\sum}_{1\leq j_1<\cdots<j_p\leq n}}_{j_1,\cdots,j_p\notin\{i_1,\cdots,i_k\}}\sum_{q=1}^k\left(\mathop{\sum\cdots\sum}_{1\leq u_1<\cdots<u_q\leq k}\lambda_{i_{u_1},\cdots,i_{u_q},j_1,\cdots,j_p}\right)\right)$$ we have $$\bar{F}_T(t)=\sum_{k=1}^n(-1)^{k-1}\mathop{\sum\cdots\sum}_{1\leq i_1<\cdots<i_k\leq n}\exp\left\{-A_{i_1,i_2,\ldots,i_k}t\right\}.$$

Hence, $$f_T(t)=-\frac{\mathrm{d}}{\mathrm{dt}}\bar{F}_T(t)=\sum_{k=1}^n(-1)^{k-1}\mathop{\sum\cdots\sum}_{1\leq i_1<\cdots<i_k\leq n}A_{i_1,i_2,\ldots,i_k}\exp\left\{-A_{i_1,i_2,\ldots,i_k}t\right\}$$

Therefore, $$E^{\bar{F}}(t)=\frac{\sum_{k=1}^n(-1)^{k-1}\mathop{\sum\cdots\sum}_{1\leq i_1<\cdots<i_k\leq n}\left[\exp\left\{-tA_{i_1,i_2,\ldots,i_k}\right\}-\exp\left\{-t\left(\sum_{p=1}^k\lambda_{i_p}\right)\right\}\right]}{1-\prod_{i=1}^n (1-e^{-\lambda_it})}$$
\subsubsection{Multivariate Gumbel's Type 1 ($MG1$)}
The resulting joint distribution of lifetimes $X_1,\ldots,X_k$ of the components having $MG1$ is given by
\begin{eqnarray}
\overline{F}(x_1,x_2,\ldots,x_n)&=&\exp\Big\{-\sum_{i=1}^{n}\lambda_i t_i-\sum_{i_1<}\sum_{i_2}\lambda_{i_1,i_2}x_{i_1}x_{i_2}-\sum_{i_1<}\sum_{i_2<}\sum_{i_3}\lambda_{i_1,i_2,i_3}x_{i_1}x_{i_2}x_{i_3}-\ldots-\nonumber\\
&&\lambda_{i_1,i_2,\ldots,i_n}x_{i_1}x_{i_2}\ldots x_{i_n}
\Big\}\nonumber.
\end{eqnarray}
\begin{equation*}
	\begin{split}
		\bar{F}_T(t)=&\sum_{i=1}^n\exp\left\{-t\lambda_i\right\}\\
		&-\mathop{\sum\sum}_{1\leq i_1<i_2\leq n}\exp\left\{-t(\lambda_{i_1}+\lambda_{i_2})-t^2\lambda_{i_1,i_2}\right\}\\
		&+\mathop{\sum\sum\sum}_{1\leq i_1<i_2<i_3\leq n}\exp\left\{-t(\lambda_{i_1}+\lambda_{i_2}+\lambda_{i_3})-t^2(\lambda_{i_1,i_2}+\lambda_{i_1,i_3}+\lambda_{i_2,i_3})-t^3\lambda_{i_1,i_2,i_3}\right\}\\
		&\vdots\\
		&+(-1)^{n-1}\exp\left\{-t\sum_{i=1}^{n}\lambda_i-t^2\mathop{\sum\sum}_{1\leq i_1<i_2\leq n}\lambda_{i_1,i_2}-t^3\mathop{\sum\sum\sum}_{1\leq i_1<i_2<i_3\leq n}\lambda_{i_1,i_2,i_3}-\cdots-t^n\lambda_{i_1,i_2,\ldots,i_n}\right\}\\
		=&\sum_{k=1}^n(-1)^{k-1}\mathop{\sum\cdots\sum}_{1\leq i_1<\cdots<i_k\leq n}\exp\left\{-\sum_{p=1}^k t^p\mathop{\sum\cdots\sum}_{1\leq u_1<\cdots<u_p\leq k}\lambda_{i_{u_1},\cdots,i_{u_p}}\right\}
	\end{split}
\end{equation*}

Therefore, $$E^{\bar{F}}(t)=\frac{\sum_{k=1}^n(-1)^{k-1}\mathop{\sum\cdots\sum}_{1\leq i_1<\cdots<i_k\leq n}\left[\exp\left\{-\sum_{p=1}^k t^pB_{i_1,i_2,\ldots,i_k}\right\}-\exp\left\{-t\left(\sum_{p=1}^k\lambda_{i_p}\right)\right\}\right]}{1-\prod_{i=1}^n (1-e^{-\lambda_it})}$$
where $B_{i_1,i_2,\ldots,i_k}=\mathop{\sum\cdots\sum}_{1\leq u_1<\cdots<u_p\leq k}\lambda_{i_{u_1},\cdots,i_{u_p}}$

\subsection{Parallel system with components having Multivariate Weibull distribution}
\subsubsection{Independent Weibull}
The survival function of $n$ independent random variables $X_1,X_2,\ldots,X_n$ having  Weibull distribution is given by
\begin{equation}
\label{indwe}
\overline{F}(x_1,x_2,\ldots,x_n)=\exp\Big\{-\sum_{i=1}^{n}\lambda_i t_i^{\alpha_i}\Big\}.\nonumber
\end{equation}
Hence,
\begin{equation*}
    \begin{split}
        \bar{F}_T(t)&=\sum_{k=1}^n (-1)^{k-1} \mathop{\sum\cdots\sum}_{1\leq i_1<\cdots<i_k\leq n}\exp\left\{-\sum_{p=1}^k \lambda_{i_p}t^{\alpha_{i_p}}\right\},
    \end{split}
\end{equation*}
\begin{equation*}
    \begin{split}
        f_T(t)&=\sum_{k=1}^n (-1)^{k-1} \mathop{\sum\cdots\sum}_{1\leq i_1<\cdots<i_k\leq n}\exp\left\{-\sum_{p=1}^k \lambda_{i_p}t^{\alpha_{i_p}}\right\}\sum_{p=1}^k \lambda_{i_p}\alpha_{i_p}t^{\alpha_{i_p}-1}
    \end{split}
\end{equation*}

\subsubsection{Farlie-Gumbel-Morgenstern Weibull (FGMW) System of distribution}
A n-dimensional random variable $(X_{1},X_{2},\ldots,X_{n})$ is said to follow Farlie-Gumbel-Morgenstern Weibull (FGMW) distribution if its survival function is given by
\begin{equation}
\overline{F}(x_1,x_2,\ldots,x_n)=\Big\{\prod_{i=1}^{n}\overline{F}_{i}(t_i)\Big\}\Big\{1+\gamma\prod_{i=1}^{n} [1-\overline{F}_{i}(t_i)]\Big\}\nonumber
\end{equation}
where $\overline{F}_{i}(t_i)=\exp(-\lambda_i x_{i}^{\alpha_i})$ for $i=1,2,\ldots,n$ is the survival function of $X_{i}$ and $-1<\gamma<1.$ \\
The survival function of a parallel system formed out of $n$ components having FGMW are respectively given by
\begin{eqnarray}
\overline{F}^{T}_{1}(t)
&=&\sum_{i=1}^{n}\exp(-\lambda_{i}t^{\alpha_{i}})-\sum_{i_{1}<}\sum_{i_{2}}\exp\big(-\lambda_{i_{1}}t^{\alpha_{i_{1}}}-\lambda_{i_{2}}t^{\alpha_{i_{2}}}\big)+\sum_{i_{1}}\sum_{}\nonumber\\
&=&\theta(t)+(-1)^{n-1}\exp(A(t))\Big(\gamma \prod_{i=1}^{n}(1-\exp(-\lambda_i t^{\alpha_i}))\Big)=\theta(t)+(-1)^{n-1}\exp(A(t))\gamma (1-\theta (t)),\nonumber
\end{eqnarray}
where $$\theta(t)=\sum_{k=1}^{n}(-1)^{k-1}\sum_{i_{1}<} \sum_{i_{2}<} \ldots\sum_{<i_{k}}\exp\Big(\sum_{p=1}^{k}(-\lambda_{i_{p}}t^{\alpha_{i_{p}}})\Big).$$
The corresponding failure rate, reversed failure rate and aging intensity functions are given by
$$r^{T}_1(x)=\frac{-\theta^{'}(t)+(-1)^{n}\exp(A(t))\gamma (-\theta^{'}(t))+(-1)^{n}\exp(A(t))A^{'}(t)\gamma (1-\theta(t)))}{\theta(t)+(-1)^{n-1}\exp(A(t))\gamma (1-\theta(t)))}$$
$$\mu^{T}_1(t)=\frac{-A^{'}(t)\phi(t)-\phi^{'}(t)}{1-\phi(t)\exp(A(t))}, L^{T}_1(t)=t\Big\{\frac{ A^{'}(t)\phi(t)+\phi^{'}(t)}{\phi(t)\Big(A(t)+\ln \phi(t)\Big)}\Big\},$$
Here, $E^{\overline{F}}(t)=\frac{(-1)^{n-1}\exp(A(t))\gamma (1-\theta(t))}{\theta(t)}, E^{r}(t)=\frac{\theta^{'}(t)+A^{'}(t)\theta(t)(1-\theta(t))}{\theta^{'}(t)\Big\{\theta-1+\frac{(-1)^{n}\theta(t)}{\gamma\exp(A(t))}\Big\}}, E^{\mu}(t)=\frac{\theta(t)-A^{'}(t)}{\theta^{'}(t)\Big\{1+\frac{(-1)^{n}}{\gamma\exp(A(t))}\Big\}}.$
We note that  $E^{\overline{F}}(t)\geq 0$ if $\gamma$ is positive and $n$ is odd or if $\gamma$ is negative and $n$ is even. On the other hand, $E^{\overline{F}}(t)\leq 0$ if $\gamma$ is positive and $n$ is even or if $\gamma$ is negative and $n$ is odd.  $E^{r}(t)\leq 0$ if $\gamma>0$ and $n$ is odd or if $\gamma<0$ and $n$ is even. $E^{\mu}(t)\leq 0$ if $\gamma>0$ and $n$ is even or if $\gamma<0$ and $n$ is odd. $E^{\mu}(t)\geq 0$ if $\gamma>0$ and $n$ is odd or if $\gamma<0$ and $n$ is even.
\section*{Acknowledgement}
Subarna Bhattacharjee  would like to thank Odisha State Higher Education Council for providing support to carry out the research project under OURIIP, Odisha, India (Grant No. 22-SF-MT-073).

\end{document}